\documentclass[11pt,fleqn]{article}
\font \bb=msbm10 

\def\C{\hbox{\bb C}}
\def\K{\hbox{\bb K}}
\begin{document}
\begin{center}
\Large { \bf More On Generalized Symmetries
of Partial Differential Equations and Quasiexact Solvability}
\end{center}
\newtheorem{definition}{Definition}
\newtheorem{Proposition}{Proposition}
\newtheorem{Lemma}{Lemma}
\newtheorem{Theorem}{Theorem}
\begin{center}
{\it Arthur SERGHEYEV}
\end{center}
\begin{center}
{ Institute of Mathematics of the National Academy of Sciences of
Ukraine, \\
Tereshchenkivs'ka Street 3, 252004 Kyiv, Ukraine.} \\
{\it e-mail:  arthur@apmat.freenet.kiev.ua}
\end{center}
\begin{flushright}
\begin{minipage}{0.9\textwidth}
\small
For the class of systems of PDEs, for which infinitesimal
translations (with respect to some (in)dependant variables) possess
specific finite-dimensional invariant subspaces of the space of
generalized symmetries of the systemconsidered. We establish when
there exist generalized symmetries from these subspaces, which depend
from the above-mentioned variables.
\end{minipage}
\end{flushright}
\section*{Introduction}

We proceed here with the study of generalized symmetries of
systems of PDEs, started in \cite{as}. Let us remind that the basic
idea of the appoach proposed consists in considering
finite-dimensional subspaces of the space of all the generalized
symmetries and the action of infinitesimal shifts with respect to
some (in)dependant variables\footnote{we shall refer to these
variables as to "selected" ones in what follows.}
on these subspaces. This enabled us to obtain the explicit formulas,
describing the dependance of the generalized symmetries from these
subspaces on the "selected" variables. We refer the reader to
\cite{o, fn, i, s, m} for the general motivations of the study of
generalized symmetries of PDEs.

In this paper we apply the results of \cite{as} to establish a simple
criterion of existence of the generalized symmetries from the
above-mentioned finite-dimensional subspaces, which {\it really\/}
depend on the "selected" variables.

As in \cite{as}, let us consider the system of PDEs of the form:
\begin{equation}  \label{1}
F _{\nu}(x,u,\dots, u^{(d)}) = 0, \quad \nu =1, \dots ,f,
\end{equation}
where $u=u(x)=(u_1, \dots, u _n)^{T}$
is unknown vector-function of $m$
independant variables $x=(x_1, \dots, x_m)$;
$u^{(s)}$
denotes the set of derivatives of $u$ with respect to $x$ of the
order $s$; $^{T}$ denotes matrix transposition.

\begin{definition}{\rm \cite{o}}. The differential operator ${\rm Q}$
of the form
\begin{equation} \label{2}
{\rm Q} = \sum _{i=1} ^{m} \xi _i (x,u,\dots, u^{(q)}) \partial
/\partial x_i +
\sum _{\alpha =1} ^{n} \eta _{\alpha} (x,u,\dots, u^{(q)}) \partial /
\partial u_{\alpha}
\end{equation}
is called the generalized symmetry
of order $q$ of the system of PDEs
(\ref{1}) if
its prolongation ${\rm {\bf pr}}{\rm Q}$ annulates (\ref{1}) on the
set $M$ of (sufficiently smooth) solutions of (\ref{1}):
\begin{equation} \label{3}
  {\rm {\bf pr} Q} [F_ {\nu}] \mid _M = 0, \quad  \nu =1, \dots ,f.
\end{equation}
\end{definition}

Let (like in \cite{as}) $Sym$ denote Lie algebra over the field $\C$ of
complex numbers\footnote{In fact, everywhere in our considerations
(except the example) $\C$ may be replaced by the arbitrary
algebraically closed field $\K$, since Theorem 1 below holds true for
this case too \cite{as}.}
(with respect to the so-called Lie bracket $[,]$ \cite{o})
of all the generalized symmetries of (\ref{1}) of non-negative orders,
$Sym ^{(q)}$ be the linear space of the generalized symmetries of (\ref{1})
of order not higher than $q$, $Sym _{q} \equiv Sym^{(q)}/Sym^{(q-1)}$
$(q \neq 0)$, $Sym _0 \equiv Sym^{(0)}$.

Let us denote the coordinates on the manifold of 0-jets $M^{(0)}$
\cite{o} as $z_A$: $z_{i} = x_{i}, i=1, \dots, m, z_{m+\alpha } = u _{\alpha},
\alpha =1, \dots, n$ (from now on the indices $A, B, C, D, \dots$
will run from 1 to $m+n$ and we shall denote $\partial _A \equiv
\partial /\partial z_A$).

In \cite{as} we have proved the following statement:
\begin{Theorem} \label{t4}
Let $W$ be a linear subspace of the linear space of all the
differential operators of the form (\ref{2}) of arbitrary finite
orders, $A_1, \dots , A_g$ be fixed integers from the range $1,
\dots, m+n$, $g \leq m+n$,
$W_{A_{1}, \dots, A_{g}}$ be the linear space of the operators, obtained from
the operators from $W$ by setting $z_{A_{1}}=0, \dots, z_{A_g} = 0$ in their
coefficients, $V= W \bigcap Sym$, for some $q_1$ the dimension of the
subspace of $V$ $V^{(q_{\lower2pt\hbox{\tiny{1}}})} \equiv W \bigcap Sym
^{(q_{\lower2pt\hbox{\tiny{1}}})}$
$v^{(q_{\lower2pt\hbox{\tiny{1}}})} < \infty$
and for any generalized symmetry of (\ref{1}) ${\rm Q} \in
V^{(q_{\lower2pt\hbox{\tiny{1}}})}$ $\partial {\rm Q} / \partial z_{A_s}
 \in V ,s=1, \dots, g$.

Then in each $V^{(q)} \equiv W \bigcap Sym^{(q)}$ ($q=0,
\dots, q_1)$ there exists
a basis of linearly independant
generalized symmetries ${\rm Q}_l^{(q,\gamma)}$,
$l=1, \dots, r_{\gamma}^{(q)}$, $\gamma=1, \dots, \rho ^{(q)}$ ($\rho^{(q)}
\leq v^{(q)}$, $\sum_{\gamma =1}^{\rho^{(q)}} r_{\gamma}^{(q)}
= v^{(q)}$) of the form
\begin{equation} \label{8aaa}
\begin{array}{lll}
{\rm Q}_l^{(q,\gamma)} =
\exp(\sum\limits_{s=1}^{g}
\lambda _{\gamma} ^{(q,A_{s})} z_{A_{s}})
\times \\
\sum\limits_{j_{1}=0} ^{k_{\gamma}^{(q,A_{1})} - 1}
\dots \sum\limits _{j_{g}=0} ^{k_{\gamma}^{(q,A_{g})} - 1}
(z_{A_{1}}) ^{j_{1}} (z_{A_{2}}) ^{j_{2}} \dots (z_{A_{g}}) ^{j_{g}}
\: {\rm C}_{l,j_{1}, \dots, j_{g} }^{(q,\gamma)},
\end{array}
\end{equation}
where ${\rm C}_{l,j_{1}, \dots , j_{g} }^{(q,\gamma)}$
are some differential operators from $W_{A_{1},\dots, A_{g}}$
of order $q$ or lower;
$\lambda _{\gamma} ^{(q,A_{s})} \in \C$ are some constants and
$k_{\gamma}^{(q,A_{s})}, s=1, \dots, g$ are some fixed numbers from the
range $1,\dots, r_{\gamma}^{(q)}$.
\end{Theorem}

\section{The criterium of dependance from "selected" variables}

 First of all let us notice that acting on any ${\rm Q}_{l}^{(q,\gamma)}$
(\ref{8aaa})
by the operators $\partial /\partial z_{A_s} - \lambda_{\gamma}^{(q,A_s)}$,
$s=1,\dots,g$ appropriate number of times, we can obtain the generalized
symmetry from $V^{(q)}$ of the~form
\begin{equation} \label{specsym}
\begin{array}{lll}
{\rm R}_{l}^{(q,\gamma)} = \exp(\sum\limits_{s=1}^{g}
\lambda_{\gamma}^{(q,A_{s})} z_{A_{s}})
\times \\
\sum\limits_{j_{1}=0}^{\varepsilon_{\gamma}^{(q,A_1)}}
\dots \sum\limits_{j_{g}=0}^{\varepsilon_{\gamma}^{(q,A_g)}}
(z_{A_{1}}) ^{j_{1}} (z_{A_{2}}) ^{j_{2}} \dots (z_{A_{g}}) ^{j_{g}}
{\rm K}_{l, j_{1}, \dots , j_{g} }^{(q,\gamma)},
\end{array}
\end{equation}
where
${\rm K}_{l, j_{1}, \dots ,
j_{g}}^{(q,\gamma)}$ are differential operators from $W_{A_{1},
\dots, A_{g}}$ of order $q$ or
lower, $\varepsilon_{\gamma}^{(q,A_s)} = 0$ if
$\lambda_{\gamma}^{(q,A_s)} \neq 0$ and
$\varepsilon_{\gamma}^{(q,A_s)} = 1$ if $\lambda_{\gamma}^{(q,A_s)} =
0$, $s=1, \dots,g$.

 Moreover, if
${\rm Q}_{l}^{(q,\gamma)}$ (\ref{8aaa})
really depends from $z_{A_i}$, acting on it $k_{\gamma}^{(q,A_{i})} -
1-\varepsilon_{\gamma} ^{(q,A_i)}$ times
by the operator $\partial /\partial z_{A_i} - \lambda_{\gamma}^{(q,A_i)}$ and
appropriate number of times by the operators $\partial /\partial
z_{A_s} - \lambda_{\gamma}^{(q,A_s)}$, $s \neq i$, yields the
generalized symmetry of the form (\ref{specsym}), where at least one
operator ${\rm K}_{l, j_{1}, \dots , j_{g} }^{(q,\gamma)}$ with
$j_i=\varepsilon_{\gamma}^{(q,A_i)}$ is not equal to zero.

Let us mention that the above procedure of obtaining the symmetries
of the form (\ref{specsym}) from (\ref{8aaa}) gives not the unique
symmetry (\ref{specsym}) but the set of such symmetries. We will show
how to obtain one of such symmetries, satisfying the conditions,
imposed on ${\rm K}_{j_1,\dots, j_g}^{(q,\gamma)}$ above. Let (for
the sake of simplicity) $i=1$.
Among the terms in (\ref{8aaa}) with $j_1=k_{\gamma}^{(q,A_1)} - 1 \equiv
r_1$ we choose the term(s) with maximal value of $j_2$, which we'll denote
$r_2$; among the last ones we choose the term(s) with maximal value
of $j_3$, which will be
denoted as $r_3$ and so on. At the end of this procedure we will
obtain the only (nonzero!) term from (\ref{8aaa}), in which
$j_s=r_s$, $s=1, \dots, g$. Then, if we act on ${\rm
Q}_l^{(q,\gamma)}$ (\ref{8aaa}) by the operator
$$
\prod\limits_{s=1}^{g} (\partial /\partial z_{A_s} -\lambda_{\gamma}
^{(q, A_s)})^{r_{s} - \varepsilon_{\gamma}^{(q,A_{s})}},
$$
 the presence of this term will assure that in resulting
expression (which obviously will be of the form (\ref{specsym})) at least one
${\rm K}_{l, j_{1}, \dots , j_{g} }^{(q,\gamma)}$ with
$j_1=\varepsilon_{\gamma}^{(q,A_1)}$ is not equal to zero.
Similar considerations may be undertaken for all other values of $i$.

The possible non-uniqueness of this construction follows from the fact that
(we again restrict to $i=1$) it is possible to act in a different
order: e.g. instead of searching the terms with maximal $j_2$ to
search the terms with maximal $j_g$, amongst them -- those with
maximal $j_3$, and only at the end consider $j_2$.

Now let us formulate explicitly the statement to be proved:
\begin{Theorem} \label{t5}
Provided the conditions of Theorem \ref{t4} are fulfilled, the system
(\ref{1}) possesses generalized symmetries from $V^{(q)}$, $q=0,\dots, q_1$
with $z_{A_i}$-dependant coefficients if and only if there exists
generalized symmetry of (\ref{1}) from $V^{(q)}$ of the form
\begin{equation} \label{specsym2}
{\rm R} = \exp(\sum\limits_{s=1}^{g}
\lambda_s z_{A_{s}})
\sum\limits_{j_{1}=0}^{\varepsilon_1}
\dots \sum\limits_{j_{g}=0}^{\varepsilon_g}
(z_{A_{1}}) ^{j_{1}} (z_{A_{2}}) ^{j_{2}} \dots (z_{A_{g}}) ^{j_{g}}
{\rm K}_{j_{1}, \dots , j_{g} },
\end{equation}
where $\lambda_s \in \C$,  $\varepsilon_s =
0$ if $\lambda_s \neq 0$ and $\varepsilon_s = 1$ if $\lambda_s = 0$,
$s=1, \dots,g$ and ${\rm K}_{j_{1}, \dots,j_{g}}$ are
differential operators from $W_{A_{1},\dots, A_{g}}$ of order $q$ or
lower, and for this symmetry either $\lambda _i \neq 0$ or $\lambda
_i=0$ but at least one operator ${\rm K}_{j_{1}, \dots , j_{g}}$ with
$j_i=1$ is not equal to zero.
\end{Theorem}
{\it Proof.}
Sufficiency is obvious. Necessity follows from the fact that the
basis of linearly independant generalized symmetries of (\ref{1}) from all
$V^{(q)}$, $q=0, \dots, q_1$ is given by (\ref{8aaa}) in virtue
of Theorem \ref{t4} and from the above considerations. $\triangleright$

Thus, in order to check the existence of $z_{A_i}$-dependant
symmetries of (\ref{1}) from $V^{(q_{\lower2pt\hbox{\tiny{1}}})}$ it
suffices (provided the conditions of Theorem \ref{t4} are fulfilled)
to check the existence of the specific symmetries (\ref{specsym}), as
described above in Theorem \ref{t5}.

In particular, for $g=1$ (\ref{1}) admits $z_{A_1}$-de\-pendant
symmetries from $V^{(q)}$, $q=0, \dots, q_1$ if and only if there
exists the symmetry from $V^{(q)}$ of the form
\begin{equation} \label{cases}
{\rm Q}= \exp(\lambda z_{A_1}) {\rm K}_{0}, \lambda \in \C, \lambda \neq 0
\: \mbox{or} \:
{\rm Q}= {\rm K}_{0} + z_{A_1} {\rm K}_{1}, {\rm K}_1 \neq 0,
\end{equation}
where ${\rm K}_{0}$ and ${\rm K}_{1}$ are the operators from $W$ with
$z_{A_1}$-independant coefficients of order not higher than $q$.

Now let us illustrate all these ideas by the following

{\it \bfseries Example.} Let $m=2$, $n=1$, $u_1 \equiv u$, $x=(x_1 \equiv t,
x_2 \equiv y)$, $u _{(l)} = \partial ^{l} u / \partial y ^{l}$ and
(\ref{1}) be the evolution equation
\begin{equation} \label{eveq}
\partial u / \partial t =G(u, u_{(1)}, \dots, u _{(d)}), \quad d\geq 2,
\end{equation}
$W$ be the linear space of the differential operators of the form (\ref{2})
with $\xi _i \equiv 0$, whose coefficient
$\eta \equiv \eta _{1}$, which is called the characteristic of the symmetry
\cite{s},
depends only on $y, u, u_{(1)}, u_{(2)}, \dots $, $V = W \bigcap Sym$.
In \cite{s} it is proved that for such a $V$ $v^{(q)} \leq v^{(1)} + q - 1$
for $q=1,2, \dots$ and $v^{(1)} \leq d+3$.

 In \cite{as} it is proved that the generalized
symmetries of order $q \geq 2$ from $W$ depend on $y$ as polynomials of
order not higher than $v^{(q)} -1 \leq v^{(1)}+q-2$.
Therefore, according to Theorem \ref{t5} and (\ref{cases}),
the equation (\ref{eveq}) posseses $y$-dependant symmetries from $W$
if and only if there exists a symmetry from $W$, which depends on $y$
linearly.

As a final remark, let us mention that all the
generalizations of Theorem \ref{t4}, presented in \cite{as}, imply
the corresponding (evident) generalizations of Theorem \ref{t5}.

\end{document}